\documentclass[10pt]{amsart}
\usepackage{amsmath,amssymb,amsthm,amscd,amsopn,amsfonts,amsxtra}
\usepackage{latexsym,array,verbatim}
\usepackage[mathscr]{eucal}
\usepackage{eucal} 
\usepackage{eufrak} 
\usepackage{graphicx,color}
\usepackage{epsfig}
\numberwithin{equation}{section}
\newtheorem{theorem}{Theorem}[section]
\newtheorem{lemma}[theorem]{Lemma}
\newtheorem{corollary}[theorem]{Corollary}
\newtheorem{proposition}[theorem]{Proposition}
\newtheorem{assumption}[theorem]{Assumption}

\theoremstyle{definition}
\newtheorem{problem}[theorem]{Problem}

\definecolor{orange}{rgb}{0.995, 0.75, 0.35}
\definecolor{purple}{rgb}{0.7, 0.2, 0.5}
\definecolor{royalblue}{rgb}{0.2, 0.7, 0.8}
\definecolor{darkgreen}{rgb}{0.2,0.725,0.25} 

\def\al{\alpha}
\def\de{\delta}
\def\eps{\epsilon}
\def\ga{\gamma}
\def\lam{\lambda}

\def\vphi{\varphi}
\def\De{\Delta}
\def\Ga{\Gamma}

\def\iy{\infty}
\def\pa{\partial}

\def\inv{^{-1}}

\def\sh{\sinh}

\def\sech{\mathrm{sech}}

\def\supp{\mathrm{supp}}

\def\H{\sqrt{H}}

\newcommand{\la}{\langle}
\newcommand{\ra}{\rangle}

\newcommand{\hr}{\hookrightarrow}
\newcommand{\hB}{\hfill$\Box$}
\newcommand{\mcal}{\mathcal}

\newcommand{\rk}{{\bf Remark.}\ \ }

\newcommand{\nd}{\noindent}
\newcommand{\n}{\newline}

\newcommand{\vs}{\vspace}

\newcommand{\Z}{\mathbb{Z}}
\newcommand{\R}{\mathbb{R}}
\newcommand{\C}{\mathbb{C}}
\newcommand{\N}{\mathbb{N}}

\def\sideremark#1{\ifvmode\leavevmode\fi\vadjust{\vbox to0pt{\vss
 \hbox to 0pt{\hskip\hsize\hskip1em
\vbox{\hsize2cm\tiny\raggedright\pretolerance10000
 \noindent #1\hfill}\hss}\vbox to8pt{\vfil}\vss}}}%

{\end{list}}

\begin{document}
\title[Harmonic analysis related to Schr\"odinger operators]
{Harmonic analysis related to Schr\"odinger operators}
\author{Gestur Olafsson}
\address[Gestur \'Olafsson]{Department of Mathematics \\
         Louisiana State University  \\
         Baton Rouge, LA 70803}
\email{olafsson@math.lsu.edu}      
\urladdr
{http://www.math.lsu.edu/\textasciitilde{olafsson}}
\author{Shijun Zheng}
\address[Shijun Zheng]{Department of Mathematical Sciences\\
Georgia Southern University\\
Statesboro, GA 30460}
\address{ and}
\address{Department of Mathematics \\
          University of South Carolina  \\
         Columbia, SC 29208}
\email{szheng@georgiasouthern.edu}
\urladdr
{http://math.georgiasouthern.edu/\symbol{126}{szheng}}
\thanks{G.~Olafsson is supported by NSF grant DMS-0402068.  S.~Zheng is supported by DARPA grant HM1582-05-2-0001}         
\keywords{Besov space, Littewood-Paley theory, 
Schr\"odinger operator}
\subjclass[2000]{Primary: 42B25; Secondary: 
35J10} 
\date{October 10, 2007}
\begin{abstract}
In this article we give an overview on  
some recent development of Littlewood-Paley theory for Schr\"odinger operators. 
We extend the Littlewood-Paley theory for special potentials considered in
the authors' previous work.  
We elaborate our approach by considering 
potential in $C^\infty_0$ or Schwartz class in one dimension.  
In particular the low energy estimates are treated by establishing 
some new and refined asymptotics for the eigenfunctions and their Fourier transforms. 
We give maximal function characterization of the Besov spaces and
Triebel-Lizorkin spaces associated with $H$.  
We then prove a spectral multiplier theorem on these spaces and 
derive Strichartz estimates for the wave equation with a potential.
We also consider similar problem for the unbounded potentials in the Hermite and Laguerre cases,
whose $V=a|x|^2+b|x|^{-2}$ are known to be critical in the study of perturbation of nonlinear dispersive equations. 
This improves upon the previous results 
when we apply the upper Gaussian bound for the heat kernel and its gradient.
\end{abstract}

\maketitle 


\section{Introduction}

The purpose of this article is to review recent development of harmonic analysis for differential operators, in purticular a 
Schr\"odinger operator $H=-\De+V$, where $V$ 
is a real-valued potential function on $\R^n$.
We are interested in developing the Littlewood-Paley 
theory for $H$ in order to understand the associated 
function spaces and their roles in 
dispersive partial differential equations. 


This subject has been drawing increasing 
attention in the area of harmonic analysis and PDE 
\cite{JN94,D97, 
DOS02,DP05,DZ05, E95,E97,FJW,He90a, 
BZ,OZ06, Sch05b,GS04,BT06}, 
to name only a few.  
The function space theory for $H$ was introduced in \cite{E95,E97,D97} for the
Hermite and Laguerre operators. In \cite{BZ,OZ06} the authors considered 
Littlewood-Paley theory for $H$ with special potentials in an effort to extend the
function space theory to the bounded potential case. In this paper we will 
summerize and develop the fundamental theory for general Schr\"odinger operators on $\R^n$. 
Furthermore we obtain a  Littlewood-Paley decomposition for $L^p$ spaces as well as Sobolev spaces using dyadic functions of $H$. 
We elaborate our approach by considering 
one dimensional $H$ with $V$ in 
$\mcal{S}(\R)$, the Schwartz class.  We will give outlines of the proofs for some of the main results and refer to the 
references, either old or new, for the detailed proofs.

\subsection{Besov and Triebel-Lizorkin spaces}
For a (Borel) 
measurable function $\phi$: $\R\to\C$ we define by functional calculus 
\[ \phi(H)=\int_{-\infty}^\infty  \phi(\lam) dE_\lam\,  \]
where $H=\int \lam dE_\lam$ is the spectral resolution of $H$.

Let $\{\varphi_j\}_{-\iy}^\iy\subset C_0^\infty ({\R}) $ be a dyadic system satisfying 
\begin{enumerate} 
\item[(i)]  
 $\supp\; \varphi_j
\subset \{ x: 2^{j-2}\le |x|\le 2^j\} $,
\item[(ii)] $|\vphi_j^{(k)}(x)|\le c_k 2^{-kj}\, ,    \qquad \forall j\in \Z$, $k\in \N_0=\{0\}\cup\N$,
\item[(iii)] $\displaystyle{
\sum_{j=-\infty}^\infty |\varphi_j(x)| \approx 1, \quad \forall x\neq 0\,.}$
\end{enumerate}

\nd
Let $\alpha\in \R$,  $0<p\le \infty, 0<q\le \infty$ and
$\{\varphi_j\}_{j\in\Z} $ be as above. The homogeneous \emph{Besov space} \emph{associated with} $H$, denoted by $\dot{B}_p^{\alpha,q}(H)$, %
is defined to be the completion of $\mathcal{S}(\R^n)$ with respect to  
the quasi-norm 
\begin{align}\label{e:homo-B-H}
\Vert f\Vert_{\dot{B}_p^{\alpha,q}(H) } = 
\big(\sum_{j=-\infty}^{\infty}
 2^{j\alpha q} \Vert \vphi_j(H)f \Vert_{L^p}^q \big)^{1/q}\, .
\end{align}


Similarly, the homogeneous {\em Triebel-Lizorkin space associated with $H$}, denoted by $\dot{F}_p^{\alpha,q}(H)$, $\alpha\in \R$,  $0<p< \infty, 
0<q\le \infty$ is defined by the quasi-norm
\begin{equation*}\label{eq:F-norm}
\Vert f\Vert_{\dot{F}_p^{\alpha,q}(H)} 
=\Vert \big(\sum_{j=-\infty}^{\infty} 2^{j\alpha q} \vert \vphi_j(H)f \vert^q\big)^{1/q}\Vert_{L^p} \, .
\end{equation*} 


We are mainly concerned with the following three interrelated problems.
\begin{problem}\label{p:a-b-c}
\begin{itemize}
\item[a.] Littlewood-Paley theory for $B(H)$ and $F(H)$ 
using maximal function characterization.
\item[b.] Spectral multiplier theorem: Find 
sufficient condition for $\mu\in L^\iy(\R)$  
such that $\mu(H)$ is bounded on $L^p$, 
$B(H)$, and $F(H)$. 
\item[c.] 
Strichartz estimates for $e^{-it\sqrt{H}}$ and 
$e^{-itH}$ 
which measure the spacetime regularity of solutions to  
wave and Schr\"odinger equations.
\end{itemize} 
\end{problem}

The decay estimate in (\ref{e:der-ker-phi}) has been known to be fundamental and useful in function space theory and spectral multiplier problem \cite{E95,D97,BZ,OZ06,OOZ06}. 
We will see that it can be applied to characterize 
$B(H)$ and $F(H)$ spaces with full ragne of parameters $0<p,q<\iy$ and
show Mihlin-H\"ormander type multiplier result
on $L^p$, $B(H)$ and $F(H)$ spaces; 
for the multiplier problem we actually formulate a more general condition as in
 (\ref{e:poly-dec-zeta}).

Let $\phi(H)(x,y)$ denote the integral kernel of $\phi(H)$. 
\begin{assumption}\label{a:phi-dec}  
Let $\phi_j\in C_0^\infty(\R)$ satisfy the conditions in (i), (ii).  
Assume that for $\ell=0,1$ and 
 for every $N\in\N_0$ there exists a constant 
 $c_N>0$ such that for all $j\in\Z$
\begin{equation}\label{e:der-ker-phi}
\vert  \nabla_x^\ell \phi_j(H)(x,y)\vert 
\le c_N\frac{2^{(n+\ell)j/2}}{(1+2^{j/2}|x-y|)^N}\,. 
\end{equation}
\end{assumption}


We will outline the proof of the fact in Section \ref{s:V-1d} that on the real line $H$ satisfies Assumption \ref{a:phi-dec}
for $V\in C^\infty_0(\R)$. We discover that
$V$ being compactly supported and $H$ having no resonance at zero
are necessarily and sufficient 
for the gradient estimate ($\ell=1$) in (\ref{e:der-ker-phi}) to hold in low energy $-\iy<j<0$.
The kernel decay for $\ell=0$, $j=0$ 
was an open question in \cite[B.7]{Si82}. \footnote{For the sake of exposition we are
not trying to pursue how singular $V\in L^1(\R)$ can be. 
Indeed, for each $N$ the proof of Theorem \ref{th:de-phi-dec} shows that (\ref{e:der-ker-phi}) holds for $j\ge 0$ if 
 $ |V^{(s)}(t)|\lesssim\la t\ra^{-N-2-\eps}$, 
 $0\le s\le N$ and that
 (\ref{e:der-ker-phi}) holds for $j<0$ if $
 V\in L^1$ with compact support and $H$ has no resonance at zero.} 

Define the {Peetre type maximal function} for $H$ as:  for $j\in\Z$, $s>0$
$$
\vphi_{j,s}^*f(x) = \sup_{t\in \R^n} \frac{| \vphi_j(H)f(t)|}{
 ( 1+2^{j/2}|x- t  |)^{s}} \,.  \;\;
 $$
The following theorem gives a maximal  function characterization of the {\em homogeneous} 
spaces. 

\begin{theorem}\label{th:homog-B-F} Suppose $H$ satisfies Assumption \ref{a:phi-dec}
and $\{\vphi_j\}$ is a system satisfying (i)--(iii).
The following statements hold.

a) If $0< p\le\infty$, $0< q\le \infty$, $\al\in \R$ and  $s>n/ p$, then
\[ \Vert f\Vert_{\dot{B}_p^{\alpha,q}(H)}\approx
\Vert \{2^{j\al}\vphi^*_{j,s}(H)f\}\Vert_{\ell^q(L^p)} \;. \]

b) If  $0< p< \infty$, $0< q\le \infty$, $\al\in \R$ and $s>n/\min (p,q)$, then 
\[
 \Vert f\Vert_{\dot{F}_p^{\alpha,q}(H)}\approx 
\Vert \{2^{j\al}\vphi^*_{j,s}(H)f\}\Vert_{L^p(\ell^q)} \;.
\]
\end{theorem}

It is well-known that such a characterization implies that any two dyadic systems satisfying (i)--(iii)
 give rise to equivalent norms on $\dot{B}_p^{\alpha,q}(H)$ and $\dot{F}_p^{\alpha,q}(H)$. 
Another consequence is that $H$ has the lifting property. 
\begin{corollary}\label{c:lift}  
Suppose $H$ satisfies Assumption \ref{a:phi-dec}. 
Let $s,\al\in\R$ and $0<p,q\le\infty$. 
Then $H^s$ maps 
$\dot{B}_{p}^{\al,q}(H)$ isomorphically and continuously onto 
$\dot{B}_{p}^{\al-s,q}(H)$. Moreover,  
$\Vert H^sf\Vert_{\dot{B}_{p}^{\al-s,q}(H)} \approx
\Vert f \Vert_{\dot{B}_{p}^{\al,q}(H)} $. 
The analogous statement holds for $\dot{F}_{p}^{\al,q}(H) $. 
\end{corollary}

The proofs of Theorem \ref{th:homog-B-F} and Corollary \ref{c:lift} are quite standard and can be 
found in \cite{Z06a, OZ06}; see also \cite{E95, D97}. 
Using Calder\'on-Zygmund decomposition and Assumption 
\ref{a:phi-dec}  we can show that  $L^p(\R^n)=\dot{F}^{0,2}_p(H)$ 
 if $1<p<\infty$ and obtain the Littlewood-Paley inequality for $L^p$ spaces.
  If in addition $V\in S$, we can using lifting property of $H$ to 
characterize the 
Sobolev spaces ${W}^{2\al}_p(\R^n)={F}^{\al,2}_p(H)$, the inhomogeneous versions 
of $\dot{F}^{\al,2}_p(H)$, 
with equivalent norms \cite{Z06a,OZ06}.
\begin{theorem}\label{th:L-P-Lp} Let $1<p<\infty$.   The following statements hold.\n
a) If $H$ satisfies Assumption \ref{a:phi-dec}, then 
\begin{equation*}
\Vert f\Vert_{L^p(\R^n)} \approx 
\Vert \big(\sum^\infty_{j=-\infty} |\vphi_j(H)f(\cdot)|^2\big)^{1/2} \Vert_{L^p(\R^n)}\;.
\end{equation*}
b) If, in addition to the condition in a), $V\in S(\R^n)$, then for all $\al\in\R$ 
\begin{equation*}
 \Vert f\Vert_{{W}_p^{2\al}(\R^n)} \approx 
\Vert \Phi(H)f\Vert_{L^p(\R^n)}+ \Vert \big(\sum^\infty_{j=1} 2^{2j\al}|\vphi_j(H)f(\cdot)|^2\big)^{1/2} \Vert_{L^p(\R^n)}\,,
\end{equation*}
where $\{\Phi,\varphi_j\}_{j=1}^\iy$ is an inhomogeneous system satisfying (i), (ii) and ($iii^\prime$) below.
\end{theorem}
\footnote{ Again we are not trying to seek optimal condition on $V$; one can show that if $|\partial_x^kV(x)|\le c_k$, $|k|\le 2m_0-2$ 
for some $m_0\in \N$, then b) is true for $|\al|\le m_0$.  
} 

The analogous results above also hold for the inhomogenous spaces $B^{\al,q}_p(H)$, $F^{\al,q}_p(H)$ 
if using the system $\{\Phi,\varphi_j\}_{j=1}^\iy\subset C_0^\infty ({\R}) $
with $\supp\,\Phi\subset [-1,1]$, $\vphi_j$ satisfying (i), (ii) and instead of (iii)
\begin{enumerate}
\item[(iii${}^\prime$)]  $\displaystyle{
|\Phi(x)|  +\sum_{j=1}^\infty |\varphi_j(x)| \approx 1, \quad   \forall x .}$
\end{enumerate} 
However, the homogeneous space,  
which contains both high and low energy analysis of $H$, 
are essential and more useful in proving Strichartz inequality for 
wave equations; see e.g., \cite{Sch05b,KT98} and Section \ref{s:stri-H}. 

\subsection{Spectral multipliers} 
For Problem \ref{p:a-b-c} b, 
Mihlin-H\"ormander type spectral multipliers for $H$ have been considered in
\cite{He90a,E96,D01,
DOS02} and more recently \cite{Sch05b}.  
 
  In the classical case $H_0=-\De$, one can use Calder\'on-Zygmund lemma to prove the  
  $L^p$ Fourier multiplier theorem by
  showing that the kernel of $\mu(H_0)$ 
  verifies the H\"ormander condition 
\begin{equation}\label{e:hor-con}
\int_{|x-y|>2|y-\bar{y}|} |K(x,y)-K(x,\bar{y})|dx\le C
\end{equation}
if $\mu$ satisfies certain smoothness condition; 
see e.g. \cite{Hor60,St93}.
  However, for a general elliptic operator (\ref{e:hor-con}) is not available.  
   For the Schr\"odinger operator $H$ 
   with $V\ge 0$, Hebisch \cite{He90a} 
used heat kernel estimates (h.k.e) to prove a spectral multiplier theorem on $L^p$. Later on the heat kernel approach  
has been further developed  to deal with positive selfadjoint differential operators \cite{DOS02}. 
  
The question remains if the negative part of $V$ is nonzero, in which case
the upper Gaussian bound for $e^{-tH}$ may  {\em not} be valid.  
In \cite{Z06b} we are able to treat general $V$   
by replacing the h.k.e. with a (much) 
weaker condition for the pointwise decay of a spectral kernel, namely
(\ref{e:poly-dec-zeta}).

As in \cite{Z06b} the following hypothesis on $H$ 
is the main ingredients in proving spectral multiplier 
theorem. 
Let 
 $\phi_j(x)=\phi(2^{-j}x)$. 
\begin{assumption}\label{a:wei-phi-dec}  
a. (Weighted $L^2$ estimate) There exists $s>n/2$ such that
\begin{equation}\label{e:s-wei-ineq}
 \sup_y \Vert | x-y|^s \phi_j(H)(x,y)\Vert_{L^2_x} 
\le c 
2^{(n/2-s)j/2}
\qquad \forall j\in\Z ,
\end{equation}  
where $c=c(\Vert \phi\Vert_{X_*^s})$. 

b. (Weighted $L^\iy$ estimate) There exists a finite measure $\zeta$ such that 
for all $j\in\Z$ 

\begin{equation}\label{e:poly-dec-zeta}
|\phi_j(H)(x,y)|\le c' \int_{\R^n} 2^{jn/2} (1+ 2^{j/2}|x-y-u|)^{-n-\eps} d\zeta(u) , 
\end{equation}
where $c'=c'(\Vert \phi\Vert_{W_{2}^{n+\eps}})$. 
\end{assumption}
 Here $X_*^{s}=\{ f\in X^{s}(\R_*): \Vert f\Vert_{X_*^{s}}=: 
\sup_{t>0}  \Vert f(t\cdot)\eta\Vert_{X^s}<\iy\}$, 
where $\eta$ is a fixed function in $C^\iy_0$ with support away from $0$, 
$\R_*=\R\setminus \{0\}$ and $X^s$ is either 
 $W_2^s(\R)$ or 
$C^s(\R)$, the H\"older class \cite{Tr83,D01}. 

In one dimension Assumption \ref{a:wei-phi-dec} is satisfied for $X^s=C^{1}$
if $V\in L^1_2(\R)$ or $V\in L^1_1(\R)$ and $H$ has no resonance
at  zero, where  $L^1_\ga:=\{V:  (1+|x|)^\ga V\in L^1\}$. 
   In three dimensions Assumption \ref{a:wei-phi-dec} is true for $X^s=W^s_{2}$, 
if $V=V_+ -V_-$ is in the Kato class with 
small Kato norm, namly $\Vert V_\pm\Vert_K<2\pi$  
 and $H$ has no resonance at zero \cite{Z07,DP05}. 

Under Assumption \ref{a:wei-phi-dec} for $H$ and the condition $\mu\in X_*^s$ for some $s>n/2$,
 the boundedness of $\mu(H)$ on $\dot{B}^{\al,q}_{p}(H)$, $1<p<\iy$ is an immediate consequence of 
the $L^p$ result in \cite{Z06b}. 
To prove that $\mu(H)$ is bounded on 
$\dot{F}^{\al,q}_{p}(H) $ we use an $L^p(\ell^q)$ vector-valued version of the proof of 
the $L^p$ result 
by applying Calder\'on-Zygmund decomposition and the 
dyadic estimates 
(\ref{e:s-wei-ineq}) and (\ref{e:poly-dec-zeta}). 
\begin{theorem}\label{th:m(H)F} \textup{\cite{OOZ06}} Suppose $H$ satisfies Assumption 
\ref{a:wei-phi-dec}. Let $X^s=C^s$ or $W^s_2$ and $\eta$ be a fixed function in $C^\iy_0(\R_*)$. 
If there exists some $s>\frac{n}{2}$ so that  
\begin{equation*}\label{e:m-hor-con}
\sup_{t>0}  \Vert \mu(t\cdot)\eta\Vert_{X^s}<\iy, 
\end{equation*}
then $\mu(H)$ is bounded on $\dot{B}^{\al,q}_p(H)$ and $\dot{F}^{\al,q}_p(H)$,  
$1<p,q<\infty$, $\al\in \R$. 
\end{theorem}

It is easy to observe that $X=C^s$ corresponds to the usual Mihlin condition and
$X=W^s_2$ the H\"ormander condition. 
That the the exponent $\frac{n}{2}$ is sharp has been noted in 
e.g., \cite{
St93,
DOS02}. 
 Note that under the conditions in Assumption \ref{a:wei-phi-dec}, which is an alternative condition than
Assumption \ref{a:phi-dec},
 Theorem \ref{th:m(H)F} implies
\begin{align*}
 &\Vert f\Vert^\phi_{\dot{B}^{\al,q}_p(H)}\approx 
 \Vert f\Vert^\psi_{\dot{B}^{\al,q}_p(H)}\\ 
&\Vert f\Vert^\phi_{\dot{F}^{\al,q}_p(H)}\approx 
 \Vert f\Vert^\psi_{\dot{F}^{\al,q}_p(H)}
 \end{align*}
given any two system $\{ \phi_j\}_{j\in\Z}$, $\{ \psi_j\}_{j\in\Z}$, 
 which is also a corollary of Theorem \ref{th:homog-B-F} as we have mentioned.
Moreover, 
by interpolation and duality we obtain from the proof of Theorem \ref{th:m(H)F} 
 that  
\begin{align*}
 \dot{F}^{0,2}_p(H)=L^p,  \qquad 1<p<\infty
\end{align*}
(see \cite{OZ06} for the inhomogeneous case), 
which is part a) of Theorem \ref{th:L-P-Lp} while under somehow more general conditions in Assumption \ref{a:wei-phi-dec}.

\vs{.080in}
\nd {\bf Remark on Assumption \ref{a:wei-phi-dec}} \quad
Assumption \ref{a:wei-phi-dec} is {\em intrinsic} in the sense that 
it only relies on the property of $H$ and is independent of the multiplier $\mu$. 
As can be seen from the proofs in \cite{Z06b,OOZ06},  
Inequalities (\ref{e:s-wei-ineq}) and (\ref{e:poly-dec-zeta}) 
are to control higher and lower energy estimates of $\mu(H)$ respectively.  
If in (\ref{e:s-wei-ineq}) letting $\zeta=\de$, the Dirac measure,
  then we obtain the following pointwise decay 
\begin{equation*}
|\phi_j(H)(x,y)|\le   c_{n,\eps} 2^{jn/2} (1+ 2^{j/2} |x-y|)^{-n-\eps} 
\end{equation*}
which is only valid, in general, for nonnegative potentials. 
This is the reason why we call (b) a ``weighted" pointwise estimate.



The remaining of the paper is organized as follows. In Section \ref{s:stri-H}, we apply the 
interpolation properties of $B(H)$, $F(H)$ 
to obtain Strichartz estimates for $e^{-it\sqrt{H}}$.  
In Section \ref{s:outline-max} we provide the outlines of the proofs of 
Theorem \ref{th:homog-B-F} and Theorem \ref{th:L-P-Lp} under Assumption \ref{a:phi-dec}.
In Section \ref{s:V-1d} we show that for any $V$ in $\mcal{S}(\R)$,  Assumption \ref{a:phi-dec}
is verified for high energy 
and for any $V$ in $C^\iy_0(\R)$, the low energy 
estimates holds in the absence of resonance.  The proofs are based on certain new and refined estimates 
for the modified Jost functions and its Fourier transforms whose details are quite lengthy and 
will appear elsewhere. 
For unbounded potentials we consider  in Section \ref{s:hermite-laguerre} %
the analogue of Theorem \ref{th:homog-B-F} for Hermite and Laguerre operators, where $V=a|x|^2+b|x|^{-2}$, by using gradient estimates for $e^{-tH}$. 
Further, we would like to mention that the literature in the area suggests 
that it is possible to consider analogous problems 
for $H$ 
with (degenerate) magnetic potentials, cf. \cite{
Tie06,Sz06a,Z07}  
\begin{align*}
-\frac{1}{2} \sum_{j=1}^n ( \partial_{x_j}+ia_jy_j)^2+ (\partial_{y_j} -ia_jx_j)^2,\qquad a_j\in \R,
\end{align*}
by following a similar approach developed here.

\section{Strichartz estimates for $H$}\label{s:stri-H}  
 It is well known that Strichartz estimates have useful 
applications in wellposedness problem for nonlinear dispersive equations 
\cite{KT98,ST02,Sch05c}. Consider the following perturbed wave equation with a potential on $\R^{1+n}$ 
\begin{equation}\label{e:wav-V}
\begin{cases}
 u_{tt}  +Hu=F(t,x)&\, \qquad\\ 
\; u(0,x)=u_0(x),\; \; u_t(0,x)=u_1(x)&\; 
\end{cases}
\end{equation} 
whose solution is given by
\begin{align*}
u(t,x)=\cos(t\H)u_0+ \frac{\sin(t\H)}{\H}u_1+
\int_0^t \frac{\sin((t-s)\H)}{\H}F(s,\cdot)ds .
\end{align*}
When $n\ge 2$, the original Strichartz estimate for (\ref{e:wav-V}) with $V=0$ reads \cite{Str77}
\begin{equation}\label{e:ori-str}
\Vert u\Vert_{L^{\frac{2n+2}{n-1}}(\R^{n+1})}
\le c \Vert f\Vert_{\dot{W}_2^{-\frac{1}{2}}}  
\end{equation}
if $u_0=0$, $u_1=f$ and $F=0$, where $\dot{W}_2^s:=\dot{W}_2^s(\R^n)$ denotes the homogeneous Sobolev
space. 
We are interested in proving Strichartz estimates for (\ref{e:wav-V}) in the Besov space scale, assuming that
the dispersive estimate (\ref{e:u-wav-disp}) is true.   
In this section 
for convenience we will use $\dot{B}^{\al,q}_p(\sqrt{H})$ instead of $\dot{B}^{\al,q}_p(H)$, whose norm
is given by (\ref{e:homo-B-H}) with $H$ replaced by $\sqrt{H}$.
Note that $\dot{B}^{2\al,q}_p(\sqrt{H})=\dot{B}^{\al,q}_p(H)$.

\begin{assumption}\label{a:u-wav-disp}  Assume that $u(t,x)$, the solution to (\ref{e:wav-V}) with $u_0=0, u_1=f$ and $F=0$,
satisfies for all $t\neq 0$
\begin{align}\label{e:u-wav-disp}
\Vert u(t,\cdot) \Vert_{\infty} 
\le c |t|^{-(n-1)/2}\Vert f\Vert_{\dot{B}_1^{\frac{n-1}{2},1}(\sqrt{H})} \,.
\end{align}
\end{assumption}
 Dispersive estimates in (\ref{e:u-wav-disp}) 
 were obtained, for instance, in \cite{Be94,Cu00} for $V$ being smooth and in \cite{DP05} for 
$V$ in the Kato class.

The idea to treat (\ref{e:wav-V}) is to combine the arguments in \cite{KT98} and \cite{GV95}
for a free wave equation.  The Littlewood-Paley decomposition seems efficient in dealing with
this type of estimates although we do not have available the scaling invariance for $H$ or $\vphi_j(H)$,
as is important and crucial in the classical case. 

For this purpose we will need some interpolation, duality and embedding properties for $\dot{B}_p^{\al,q}(\H)$, 
 which are analogues of $\dot{B}_p^{\al,q}(\R^n)$ and can be derived as corollaries of 
 Theorem \ref{th:m(H)F} \cite{OZ06,Tr83}. 
 Directly using the classical Besov spaces would encounter commuting problem.  

From the expression of $u(t,x)$ we see that it is essential to estimate 
$e^{\pm i t\H}f$. 
Note that (\ref{e:u-wav-disp}) is  equivalent to 
\begin{align}\label{e:j-disp}
\Vert u_j(t,\cdot) \Vert_{\infty} 
\le c |t|^{-(n-1)/2} 2^{j(n-1)/2}\Vert\vphi_j(\sqrt{H}) f\Vert_{1} \,,
\end{align}
where $u_j(\cdot)=\vphi_j(\H)u(t,\cdot)$.
Using (\ref{e:j-disp}) and $TT^*$ argument in 
\cite{GV95, KT98} we obtained in \cite{Z07} the following theorem. 
\begin{theorem}\label{th:stri-H} Let $n\ge 2$, $q,r\in [2,\iy]$ and $s\in\R$. Suppose $H$ satisfies the estimate in (\ref{e:u-wav-disp}).  Then the following estimates hold. \n
a) \begin{equation}\label{e:u-LqB-B}
\Vert e^{it\H}f\Vert_{L^{q}\dot{B}_{r}^{s,2}(\H)}
\le c\Vert f\Vert_{\dot{B}_2^{s+\sigma,2}(\H)}\,,
\end{equation}
where $q=2$ and $\sigma=\sigma(q,r)$ verifies the gap condition 
$\frac{1}{q}+\frac{n}{r}=\frac{n}{2}-\sigma$. \n
b) Let $I\subset\R$ be an interval. 
\begin{align*}
\Vert \int_0^t \frac{\sin((t-s)\H)}{\H} F(s,\cdot)ds 
\Vert_{L^{q}_I\dot{B}_{r}^{s,2}(\H)}
\le c \Vert F\Vert_{L^{q'}_I\dot{B}_{r'}^{s+2\sigma-1,2}(\H)}\,,
\end{align*}
where $q'$, $r'$ denote the usual H\"older conjugate exponents of $q$, $r$, 
$(q,r)\neq (2,\frac{2n-2}{n-3})$ and 
$(q,r)$ are wave-admissible, that is, 
$\frac{2}{q}+\frac{n-1}{r}\le \frac{n-1}{2}$. 
\end{theorem}

The estimates in Theorem \ref{th:stri-H} are analogous to those in the case of zero potential \cite{GV95}.  Under additional condition, 
e.g.,  assuming $H$ satisfies (\ref{e:der-ker-phi}),  we can release the restriction $q=2$ in (\ref{e:u-LqB-B}) 
by applying Besov embedding inequality. 
Observe that then the Besov space method yields a sharper 
estimate than (\ref{e:ori-str}) 
if in (\ref{e:u-LqB-B}) taking $q=r= (2n+2)/(n-1)$, $\sigma=1/2$, $s=0$ and substituting $H^{-\frac12}f$ for $f$, 
 noticing that $\dot{B}^{0,2}_{r}(\H)\hr\dot{F}^{0,2}_{r}(\H)=L^r$ provided $r\ge 2$. 


\vs{.05710in}\nd
{\bf Problem}: Does the endpoint estimate in Theorem \ref{th:stri-H} hold with $(q,r)=(2,\frac{2n-2}{n-3})$ for $n\ge 4$? 

The endpoint estimates involving $L^q_tL^r_x$-norm were proved in \cite{KT98} in dimensions $\ge 4$
(In dimensions 2 and 3 the endpoint estimates fail). 
We do not know the answer to the question in the problem
for $L^q\dot{B}^{s,2}_r$. 
Such a result would not only be sharper but technically might involve bilinear
estimates that provide insight and deeper understanding of the 
non-scaling-invariant case when $V\neq 0$. 
 We can also formulate similar result and problem for the Schr\"odinger equation with a potential;
here we would rather refer to \cite{GS04,Sch05c,OZ06,Z07} for further discussions.

\section{Outline of Proofs}\label{s:outline-max}

\subsection{Proof of Theorem \ref{th:homog-B-F} } 
As in \cite{OZ06,E95} or \cite{Tr83}, 
Theorem \ref{th:homog-B-F} is a consequence of Peetre type maximal equality (Lemma \ref{l:phi*-M}) 
and  the 
well-known $L^p(\ell^q)$-valued Fefferman-Stein maximal inequality.
The proof of Lemma \ref{l:phi*-M} 
is standard and follows from 
 Bernstein  type inequality (Lemma \ref{l:bern-der-phi}). 
  Let \[\vphi_{j,s}^{**} f (x) = \sup_{t\in \R^n} \frac{| (\nabla_t\vphi_j(H)f)(t)|}{( 1+2^{j/2}|x -t  |)^{s}} \, .\]

\begin{lemma}\label{l:bern-der-phi} 
For $s>0$, there exists a constant $c_{n,s}>0$ such  that
 for all $j\in\Z$ 
\[
\vphi_{j,s}^{**}f(x) \leq c_{n,s} 2^{j/2} 
\vphi_{j,s}^{*} f( x), \qquad \quad \; \forall f\in \mathcal{S}(\R^n).
\]
\end{lemma}

Similar to \cite{OZ06,Tr83}  Lemma \ref{l:bern-der-phi} can be easily proved using (\ref{e:der-ker-phi}) with 
$N> n+s$. 

Let $M$ denote the Hardy-Littlewood maximal function
\begin{equation}\label{e:F-S-max}
 M f (x)=  \sup_{B\ni x}\frac{1}{|B|} \int_{B}  |f(y)| \, d y
\end{equation}
where the supreme is taken over all 
 balls $B$ in $\R^n$ centered at $x$. 

\begin{lemma}\label{l:phi*-M}   
Let $0<r<\infty$ and $s=n/r$. 
 Then 
for all $j\in\Z$
\begin{equation}\label{e:phi*-M} 
 \vphi_{j,s}^*f(x) \le c_{n,r} [M({|\vphi_{j} (H)f|}^r)]^{1/r}(x),\;\qquad \; \forall f\in \mathcal{S}(\R^n) .
\end{equation}
\end{lemma}

Thus we have seen that 
the proof of Theorem \ref{th:homog-B-F} relies on the decay estimates in Assumption \ref{a:phi-dec}. 
 In Section \ref{s:V-1d} and Section \ref{s:hermite-laguerre} we will prove such estimates for $H=-\De+V$ 
with $V\in \mcal{S}$ on $\R$, 
and $V=|x|^2$ and $V=|x|^2+b|x|^{-2}$, $b\ge 0$ on $\R^n$, where for the latter, 
the Laguerre operator, initially defined on $\R_+^n$, can be regarded as
an operator acting on $L^2(\R^n)$ for ``even'' functions.




\subsection{Proof of Theorem \ref{th:L-P-Lp}} The same proof of the identification  $F^{0,2}_p(H)=L^p$ 
in the inhomogeneous case \cite[Theorem 5.1]{OZ06} gives us 
\begin{equation}\label{e:f-F-Lp}
\Vert f\Vert_{\dot{F}^{0,2}_p(H)} \approx   \Vert f\Vert_{L^p} \;, \qquad \; 1<p<\infty
\end{equation}
for any system $\{\vphi_j\}_{j\in\Z}$ satisfying (i), (ii), (iii),  by applying $L^p(\ell^2)$-valued Calder\'on-Zygmund decomposition. \hB 


For $p=1$,  
Dziuba\'nski and Zienkiewicz 
recently obtained a characterization of Hardy space associated with $H$ using the heat operator.
\begin{theorem}\label{th:hardy-V} \textup{(\cite{DZ05})} 
Let $V\in L^{\frac{n}{2}+\eps}(\R^n)$, $n\ge 3$, with $V\ge 0$ being compactly supported.
Then  
\[
\left\|f\right\|_{\mcal{H}^1_V}\approx
 \left\|f\right\|_{\mcal{H}^1_{atom}}\approx \left\|wf\right\|_{H^1\left({\Bbb R}^n\right)} \;,\]
\end{theorem}\nd 
where $\mathcal{H}^1_V=\{f\in L^1: \sup_{t>0}  |e^{-tH}f(\cdot) |  \in L^1  \}$ and the weight $w$ is defined by $w(x)= \lim_{t\to\infty} \int_{\R^n} e^{-tH}(x,y)  dy$. The norm in the atomic decomposition
is defined as
\begin{align*} 
&\Vert f\Vert_{\mcal{H}^1_{atom}}:=\inf \sum_{j} |\lam_j|,
\end{align*}
where the infimum runs over all representations $f=\sum_j\lam_j a_j$,
$a_j$ being $\mcal{H}^1_V$ atoms satisfying (i) $\supp\,a\subset
B(x_0,r):=\{x:|x-x_0|<r\}$, 
(ii) $\Vert a\Vert_{\iy}\le |B(x_0,r)|\inv$, 
(iii) $\int a(x)w(x)dx=0$.

\vs{.060in} 
\nd 
\rk
It would be interesting to obtain a 
norm characterization for $\mcal{H}^1_V$ 
of Littlewood-Paley type in the sense of Theorem \ref{th:L-P-Lp}. 
Note that if $V\neq 0$,  
then $\mcal{H}^1_{V}\neq H^1(\R^n)$, $n\ge 3$. 
Also, in the 1D case and unbounded potential case (e.g.
$V$ being a (nonnegative) polynomial),  the $p$-atoms ($p=1$)
for $\mcal{H}_V^1$ ($V\neq 0$) are not variant of the local atoms, cf. \cite{DZ05,DZ02}.


\section{$H=-\frac{d^2}{dx^2}+V$, $V$ in $\mcal{S}(\R)$ }\label{s:V-1d} 

Associated to $H$ there exists a decomposition  $L^2(\R)=\mcal{H}_{ac}\oplus \mcal{H}_{pp}$, where 
$\mcal{H}_{ac}$ is the absolute continuous subspace and $\mcal{H}_{pp}$ the pure point 
subspace of $L^2(\R)$. Let $E_{ac}$, $E_{pp}$ be the corresponding
orthogonal projections.   If $\sigma_{ac}(H)$ denotes the absolute continuous spectrum and
$\sigma_{pp}(H)$ the pure point spectrum of $H$, then $\sigma_{ac}(H)=[0,\iy)$
and $\sigma_{pp}(H)=\{-\lam_k^2\}$ is 
a finite set of eigenvalues of $H$ in $(-\iy,0)$.

\subsection{Decay estimates of $\phi_j(H)$}
For $\lam\in \R$ let $e(x,\lam)=(1+R_0(\lam^2+i0)V)\inv e^{i\lam x}$ be the Lippman-Schwinger scattering eigenfunction
and $e_k(x)$ the $L^2$ eigenfunction of $H$ with eigenvalue $-\lam_k^2$,
where $R_0(z)=(H_0-z)\inv$ is the resolvent of $H_0=-{d^2}/{dx^2}$ with 
$z\in \C\setminus [0,\iy)$. 
Then if $\phi\in C_0$, a continuous function with compact support, we have 
\begin{equation*}\label{eq:phi(H)f}
\phi(H)f(x)=\int K(x,y) f(y) dy
\end{equation*}
where $K=K_{ac}+K_{pp}$, 
\begin{equation}\label{e:ker-ac-e}
K_{ac}(x,y)= (2\pi)^{-1} \int \phi(\lam^2)e(x,\lam)\bar{e}(y,\lam) d\lam
\end{equation}
 is the kernel of $\phi(H_{ac})$, $H_{ac}=HE_{ac}$ and $K_{pp}(x,y)=\sum_k\phi(-\lam_k^2)e_k(x)e_k(y)$
 is the kernel of $\phi(H)E_{pp}$\,.   Let $u_1$, $u_2$ be two linear independent solutions of 
  (\ref{e:jost-eq}).
 $H$ is said to have resonance at zero provided that the Wronskian of $u_1,u_2$ vanishes at zero.

\begin{theorem}\label{th:de-phi-dec}  Let $\{\phi_j\}\subset C^\infty_0(\R)$ satisfy (i),(ii).  \n 
a) If $V\in \mcal{S}$, then for each $\ell\ge 0$ and $N\ge 0$ 
there exists a constant $c_{N,\ell}$ so that for all $j\ge 0$
\begin{equation}\label{e:de-phi-j-dec}
\vert\partial_x^\ell\phi_j(H)(x,y)\vert\le c_{N,\ell}  
2^{(\ell+1)j/2}(1+ 2^{j/2}|x-y|)^{-N}\, .
\end{equation}
b) If $V\in C^\iy_0$ and $H$ has no resonance at zero, then 
 (\ref{e:de-phi-j-dec}) holds for each $\ell=0,1$, $N\ge 0$ and all $-\iy<j<0$ .
\end{theorem}

Denote by $K_j$ the kernel of $\phi_j(H)$ and write $K_j=K_{j,ac}+K_{j,pp}$.
Since 
$\sigma_{pp}$ is finite and according to e.g. \cite[Theorem C.3.4]{Si82}, 
eigenfunctions of $H$ belonging to 
 $\cap_{m,N=0}^\iy \la x\ra^{-N}W^{2m}_2(\R)=\mcal{S}(\R)$,
$K_{pp}$ satisfies (\ref{e:de-phi-j-dec}) trivially.
Hence it is sufficient to deal with $K_{j,ac}$ .

Let $R_V(z)=(H-z)\inv$. 
Let $W(\lam)$ be the Wronskian of $f_+$, $f_-$, then
 for $\lam\neq 0$ 
\begin{align*} 
& R_V(\lam^2\pm i0)(x,y)=\begin{cases}
\frac{f_+(x,\pm\lam)f_-(y,\pm\lam)}{W(\pm\lam)}& x>y\\
\frac{f_+(y,\pm\lam)f_-(x,\pm\lam)}{W(\pm\lam)}& x<y\\
\end{cases}
\end{align*}
where $f_\pm(x,z)$ are the Jost functions that solve
for $\Im z\ge 0$
\begin{equation}\label{e:jost-eq}
-f''_\pm(x,z)+V(x)f_\pm(x,z)=z^2f_\pm(x,z)
\end{equation}
and satisfy
\begin{align*} 
f_\pm(x,z)\to 
\begin{cases}e^{\pm izx} &x\to\pm \iy\\ 
\frac{1}{t(z)}e^{\pm izx}+\frac{r_\mp(z)}{t(z)}e^{\mp izx} &x\to\mp \iy
\end{cases}
\end{align*}
where $t(z), r_\pm(z)$ are called the transmission and reflection coefficients 
\cite{DT79}.

Let $m_\pm(x,z)=e^{\mp iz x}f_\pm(x,z )$ be the modified Jost function. 
We obtain, from the resolvent formula of the spectral measure of $H_{ac}$,  that if $x>y$
\begin{align}\label{e:ker-phi-m} 
& \phi(H_{ac})(x,y)=\frac{1}{2\pi }\int_{-\iy}^\iy \phi(\lam^2)m_+(x,\lam)m_-(y,\lam)
t(\lam)e^{i\lam(x-y)} d\lam 
\end{align}
where $t(\lam)= -2i\lam/W(\lam)$; see e.g. \cite{GS04}.
One can prove that the above formula and (\ref{e:ker-ac-e}) coincide. 
As a result the restriction $x>y$ can be dropped. 

To analyze the kernel we need to study the analytic asymptotics of 
$m_\pm(x,k)$,   
 namely the higher order mixed 
 differentiability on $m_\pm$
concerning the smoothness and decay in $x$ and $k$, the space and
spectral parameters.
The asymptotics of $m_\pm$ and 
 their derivatives of first order 
were originally studied in \cite{DT79}. 
However our estimates and certain formulas, including those of $t(k)$ and 
the Fourier transforms of $m_\pm$, 
are more refined and delicate. 

\subsection{Derivatives of $m(x,k)$} 

\begin{lemma}\label{l:m'-dec-k} Let $V\in \mcal{S}$ and $\ell,n\ge 0$.
 Then $m\in C^\iy(\R\times\R)$ and  for all $x, k\in\R$
\begin{align*} 
&\vert\partial_x^\ell\pa_k^n m_\pm(x,k)\vert
\le c_{n,\ell}\begin{cases}
(1+\max(0,\mp x))^n& \ell\ge 1,\\
(1+\max(0,\mp x))^{n+1} &     \ell=0. 
\end{cases}
\end{align*}
If $\ell\ge 0$, $n\ge 1$ or $\ell\ge 1$, $n\ge 0$, then
\begin{align*}
&\vert\partial^\ell_x\partial^n_km_\pm(x,k)\vert\le c_{n,\ell}\frac{(1+\max(0,\mp x))^n}{|k|^n}\,.
\end{align*}
\end{lemma} 

The proof of the lemma are based on the integral equation \cite{DT79} 
\begin{align}\label{e:m-int}
m_+(x,k)=1+\int_x^{\iy} h(t-x,k)V(t)m_+(t,k)dt ,
\end{align}
where $h(t-x,k)=\int_0^{t-x} e^{2iku}du$, and
the equation for its mixed partial derivatives: 
for $ \ell=1,2,3,\dots$, $n\ge 0$
\begin{align}\label{e:form-dx-dk-m} 
\partial_x^\ell\partial_k^n m_+(x,k)=-\sum_{j=0}^n \binom{n}{j}
(2i)^j\int_x^\infty e^{2ik(t-x)}(t-x)^j\partial_t^{\ell-1}(V m_+^{(n-j)})dt ,
\end{align}
where $m_+^{(i)}=\partial_k^{i}m_+(t,k)$.

The first estimate in the lemma is also a consequence of Lemma \ref{l:B-L1} concerning the weighted $L^1$ bound 
for $B_\pm(x,y)$, the Fourier transforms of $m_\pm-1$, 
which are especially needed for low energy estimates. The Marchenco functions $B_\pm(x,y)$ are
related to $m_\pm$ via
\[ m_\pm(x,k)=1\pm\int_0^{\pm\iy} B_\pm(x,y)e^{2iky}dy\,.
\]

 

\subsection{Analycity and asymptotics of $t(k)$ and $r_\pm(k)$} The kernel of $\phi(H_{ac})$ given by 
(\ref{e:ker-phi-m}) also requires  estimates of the coefficients $t(k)$ and $r_\pm(k)$.
\begin{lemma}\label{l:Lder-j-t} 
Let $V\in \mcal{S}$. 
Then $t(k)$, $r_\pm(k)\in C^\iy$.
If $|k|\ge 1$, then
\begin{align} 
\frac{d^n}{dt^n}t(k)=
\begin{cases}
1+O(k\inv)&n=0\\
O(k^{-n-1})&n\ge 1 
\end{cases}
\end{align}
and
\begin{align} 
\frac{d^n}{dt^n}r_\pm(k)=
O(k^{-n-1}) \qquad  n\ge 0.
\end{align}
\end{lemma}

For the proof of high energy 
 $|k|\ge 1$, we use \cite[p.145]{DT79}
\begin{align} 
&t(k)\inv=1-\frac{1}{2ik}\int_{-\infty}^\infty V(t)m_+(t,k)dt\notag\\
&(r_\pm(k)+1)t(k)\inv=1\mp\int_{-\infty}^\infty h(t,\mp k)V(t)m_\mp(t,k)dt. \label{e:r-t-hVm}
\end{align}
For low energy to show $t\in C^\iy$ the formulas we use are
more delicate. 
\begin{lemma}\label{l:t-r-resonance} Let $\nu:=W(0)$.  The following formulas hold.\n
a) If  $\nu\neq 0$, then 
\begin{align*}
&t(k)= \frac{-2ik}{W(k)} \,,
\end{align*}
where $W(k)= { -2ik+{\nu}+ \int_{-\infty}^\infty V(t) dt\int_0^\infty B_+(t,y)(e^{2iky}-1) dy}$.

\nd
b) If $\nu=0$, then 
\begin{equation*}
t(k)^{-1} 
=  1-\int_{-\infty}^\infty V(t)dt \int_0^\infty \big(\int_\xi^\infty B_+(t,\eta)d\eta\big)  e^{2ik\xi} d\xi \,.
\end{equation*}
\end{lemma}

Based on Lemma \ref{l:t-r-resonance} we give the proof of the statement $t\in C^\iy$ as in Lemma \ref{l:Lder-j-t}. 
The statement $r\in C^\iy$ follows from the fact that $t\in C^\iy$ and (\ref{e:r-t-hVm}). We divide the discussions into two cases.

\nd
{\em Case (a)} $\nu\neq 0$. 
From the formula in Lemma \ref{l:t-r-resonance} (a) it is easy to see that $W(k)$ is 
$C^s$ 
if $\int_0^\infty y^s|B(t,y)|dy\le c\la t\ra^{s+1}$, but this is true provided
$V\in L^1_{s+1}$ according to Lemma \ref{l:B-L1}. Now since $W(0)\neq 0$, it follows that
$t\in C^s$  as long as $V\in L^1_{s+1}$.

\nd
{\em Case (b)} $\nu=0$. From Lemma \ref{l:t-r-resonance} (b) we have if $s\ge 1$
\begin{align*}
&(t^{-1})^{(s)}(k) 
= -\int_{-\infty}^\infty V(t)dt \int_0^\infty (2i\xi)^s e^{2ik\xi} d\xi\big(\int_\xi^\infty B_+(t,\eta)d\eta\big), 
\end{align*}
which is the Fourier transform of the function 
\[\xi\mapsto -\chi_{(0,\iy)}(\xi)\int_{-\infty}^\infty V(t)dt(2i\xi)^s\int_\xi^\infty B_+(t,\eta)d\eta.
\]
Observe that in view of Lemma \ref{l:B-L1} this function is in $L^1$ if $V\in L^1_{s+2}$.
Hence $t\inv$ is $C^s$. 
We conclude that 
\begin{align*}
t(k)=\frac{1}{t(k)\inv} 
\end{align*}
is also in $C^s$ whenever $V\in L^1_{s+2}$, since it is a basic fact that $\nu=0\iff |t(k)|\ge c_0>0$, $\forall k$; cf. e.g., \cite[Theorem 1]{DT79}.\hB

\subsection{
Fourier transform of modified Jost function}\label{s:mar} 


The estimates for $m_\pm(x,k)$, $t(k)$ and $r_\pm(k)$,
especially in the low energy, depend on the weighted $L^1$ inequalities for the Marchenko functions $B_\pm(x,y)$
and their derivatives. Recall that $L^1_{\ga}=\{V: \int (1+|y|)^\ga |V(y)| dy<\iy  \}$ and
$W^n_{1,\ga}=\{V: \int (1+|y|)^\ga |V^{(i)}(y)| dy<\iy, i=0,\dots,n  \}$.
\begin{lemma}\label{l:B-L1}  Let $s\in\N_0$. 
a) If $V\in L^1_{s+1}(\R)$,  then there exists $c=c(\Vert V \Vert_{L^1_{s+1}} )$ so that for all $x\in \R$ 
\begin{equation*}
\int_{-\iy}^\infty |y|^s | B_\pm(x,y) | dy\le c (1+\max(0, \mp x))^{s+1} . 
\end{equation*}
b) If $V\in W^{n-1}_{1,s+1}$, $n\ge 1$, 
then there exists $c=c(\Vert V \Vert_{W^{n-1}_{1,s+1}} )$ so that for all $x\in \R$
\begin{align*}
&\int_{-\iy}^\infty |y|^s|\pa_y^nB_\pm(x,y)|dy\le
c
(1+\max(0,(\mp x))^s\\
&\int_{-\iy}^\infty |y|^s|\pa_x^nB_\pm(x,y)|dy\le
c
(1+\max(0,(\mp x))^s .
\end{align*}
\end{lemma}

The proof exploits careful iterations of the Marchenko equations
\begin{equation*}\label{e:B+marchenko}
B_+(x,y)= \int_{x+y}^\infty V(t)dt+ \int_0^y dz \int_{t=x+y-z}^\infty V(t) B_+(t,z) dt
\end{equation*}

\begin{equation*}\label{e:B-marchenko}
B_-(x,y)= \int^{x+y}_{-\infty} V(t)dt+ \int^0_y dz \int^{t=x+y-z}_{-\infty} V(t) B_-(t,z) dt, 
\end{equation*}
where we note that $B_\pm(x,\cdot)$ are supported in $\R_\pm$.

\vs{.07in}
\nd
\rk The inequality in a) is an improvement of \cite[Lemma 3.2, Lemma 3.3]{DF06},
where the cases $s=0,1$ were proved via Gronwall's inequality.



\subsection{High and low energy estimates for $\phi_j(H)(x,y)$} 
Applying Lemmas \ref{l:m'-dec-k} and Lemma \ref{l:Lder-j-t} we can  
easily prove Theorem \ref{th:de-phi-dec} (a) in high and local energy via integration by parts;
here the high energy refers to $\phi_j$, $j\ge 0$ and local energy 
refers to $\Phi\in C^\iy_0$ with $\supp\,\Phi\subset [-1,1]$.

For low energy $j<0$, (b) of Theorem \ref{th:de-phi-dec}, 
 we observe that the condition
 $\pa_xm_\pm(x,0)=0$ for large $|x|$ together with $H$ having no resonance at zero 
 is necessarily and
sufficient for $\pa_x\phi_j(H)(x,y)$ to satisfy (\ref{e:de-phi-j-dec}). 
Therefore Theorem \ref{th:de-phi-dec} (b) is true if and only if $H$ has no resonance at zero and 
$m_\pm(x,0)\equiv constant=1$ for large $|x|$.
The latter can occur {\em only when $V$ has compact support}
in view of $\pa_x^2m+2ik\pa_xm=Vm$. 

\vs{.0802in}
\nd
\rk Let $V\in \mcal{S}$. If $-\iy<j<0$, then (\ref{e:de-phi-j-dec}) still holds for $\phi_j(H)(x,y)$ with $\ell=0$. 
However, if $H$ has resonance at zero or $V$ is not compactly supported, then 
(\ref{e:de-phi-j-dec}) fails for $\pa_x \phi_j(H)(x,y)$ in the low energy as we mentioned above. 
A counterexample can be found in \cite{OZ06} for $V=-\nu(\nu+1)\sech^2x$. 
For non-smooth potentials, in \cite{BZ} we are able to obtain 
an appropriate variant of the kernel decay (\ref{e:de-phi-j-dec}) 
with $j\in\Z$, $\ell=0,1$ for $V=c\chi_{[a,b]}(x)$, $c>0$, $\chi_E$ being a characteristic function of the set
$E$.  

\vs{.102in}
\nd
{\bf Problem} For non-smooth potentials, 
it is still open 
as to determine the class of $V$ 
for which (\ref{e:de-phi-j-dec}) remains valid  especially in the low energy 
case. 

\section{$H$ satisfying upper Gaussian bound}\label{s:hermite-laguerre}
In this section we prove the analogy of Theorem \ref{th:homog-B-F} for $H$ with unbounded 
potentials, namely, the Hermite
and Laguerre operators. We see that using the upper Gaussian bound
for $\nabla p_t$ of $H$ we can prove the decay in (\ref{e:der-ker-phi}) 
in a {\em simple} way, which generalize the results of 
\cite{E95,E97} and \cite{D97}. Previously Epperson studied the 
Hermite and Laguerre cases in one dimension using oscillatory
integral method; later Dziubanski used Heisenberg group technique to prove the
kernel decay  for the Hermite expansion in $n$-dimension 
and the Laguerre expansion of integral order in 1D. 

\cite{Z06a} showed that 
 Assumption \ref{a:phi-dec} is verified when $H$ satisfies the upper Gaussian bound 
 (\ref{e:der-etH-gb}) for its heat kernel.  
The proof is based on a weighted $L^1$ inequality  
which is a scaling version of
\cite[Lemma 8]{He90a}. 


\begin{proposition}\label{pr:etH-phi}   Let $\ell=0$, $1$. 
Suppose $V\ge 0$ and $H$ satisfies the upper Gaussian bound 
\begin{equation}\label{e:der-etH-gb}
|\nabla^\ell_x e^{-tH}(x,y)| \le c_n t^{-(n+\ell )/2} e^{-c  |x-y|^2/t}\,, \qquad \forall t>0.
\end{equation}
 If $\{\vphi_j\}_{j\in\Z}$ is  
 a dyadic system  
 satisfying (i), (ii), 
 then for each $N\ge 0$ 
\begin{align*}
\vert \nabla_x^\ell \vphi_j(H) (x,y) \vert 
\le c_N 2^{j(n+\ell )/2} (1+2^{j/2}|x-y|)^{-N} , \qquad \forall j\,.
\end{align*}
\end{proposition} 

\nd
\rk 
The long time gradient estimates for $\nabla p_t$ 
($t>1$ corresponding to low energy) 
is, in general, not valid for bounded $V$, 
not even for positive $V\in \mcal{S}(\R^n)$.  

\subsection{Hermite operator $H=-\De+|x|^2$}

To verifies Assumption \ref{a:phi-dec}   it is sufficient to show that
$H$ satisfies the upper Gaussian bound in (\ref{e:der-etH-gb}),   
according to Proposition \ref{pr:etH-phi}.

For $k\in\N_0$, let $h_k$ be the $k^{th}$ Hermite function with $\Vert h_k\Vert_{L^2(\R)}=1$
such that 
\[  (-\frac{d^2}{dx^2}+x^2) h_k = (2k+1)h_k\,.
\]
Then $\{h_k(x)\}_0^\infty$ forms a complete orthonormal basis (ONB) in $L^2(\R)$. 
Let $\Phi_k(x):= h_{k_1}\otimes \cdots \otimes h_{k_n}$, $k=(k_1,\dots,k_n)\in \N_0^n$. 
Then $\{\Phi_k\}$ is an ONB in $L^2(\R^n)$.

By Mehler's formula \cite[Ch.4]{Th93}, 
the heat kernel 
has the expression 
\begin{align*}
e^{-tH}(x,y)=& \sum_{k\in\N_0^n} e^{-t (n+2|k|)} \Phi_k(x) \Phi_k(y)\\
=& (2\pi \sh (2t) )^{-n/2} e^{-\frac{1}{2}\coth (2t) (|x|^2+|y|^2) + \mathrm{cosech (2t)} x\cdot y }
\end{align*}
for all $t>0$, $x, y\in \R^n$.

One can easily calculated (cf. \cite{Z06a}) that for $\ell=0, 1$ 
there exist constants $c'>0$, $0<c<1$ 
such that 
\[
|\nabla_x^\ell p_{t}(x,y) | \le c'
\begin{cases} 
t^{-(n+\ell)/2} e^{-c |x-y|^2/t} &\quad 0<t<1\\
e^{- n t} e^{-c |x-y|^2}  &\quad t\ge 1,
\end{cases}
\]
where $p_t(x,y):=e^{-tH}(x,y)$. Hence (\ref{e:der-etH-gb}) holds. 

\subsection{Laguerre operator 
$H=L_\al=-\frac{d^2}{dx^2}+x^2+(\al^2-\frac{1}{4})x^{-2}$, $\al\ge -1/2$ 
in $L^2(\R_+)$}
The Laguerre is basically a generalization of the Hermite. 
We will show that the heat kernel estimates (\ref{e:der-etH-gb}) are valid for 
$L_\al$ with $\al\ge 1/2$.



For $
k\in \N_0$ define the Laguerre function  
\begin{equation*}\label{e:M-al-k}
 M^\al_k(x)= \left( \frac{2\,\Ga(k+1)}{\Ga(k+\al+1)}\right)^{1/2}
e^{-x^2/2}x^{\al+1/2}L^\al_k(x^2),
\end{equation*}
where $L^\al_k(x)$ are the Laguerre polynomials.  Then 
\begin{align*}
L_\al M^\al_k(x)=(4k+2\al+2)  M^\al_k(x)
\end{align*}
and $\{ M^\al_k(x)\}_0^\infty$ is an orthonormal basis
in $L^2(\R_+)$; see e.g. \cite{Th93}. 

\begin{lemma} 
Let $\al\ge 1/2$. 
Then the heat kernel of $L_\al$ satisfies 
for $\ell=0,1$ there exists constants $c'>0$, $0<c<1$ such that with $n=1$
\begin{align}\label{e:etL-gb}
&\pa_x^\ell e^{-tL_\al}(x,y)\le c'
\begin{cases}
t^{-(n+\ell)/2}e^{-c|x-y|^2/t}&0<t<1\quad\\
e^{-nt}e^{-c|x-y|^2/t}&t\ge 1. \quad \\
\end{cases}
\end{align}
\end{lemma}

\begin{proof} 
Using the following Mehler type formula:  if $0<r<1$, 
\begin{align*}
\sum_{k=0}^\infty r^k M^\al_k(x)M^\al_k(y)
=& 2e^{-i\frac{\pi}{2}\al} (xy)^{1/2} (1-r)^{-1}r^{-\al/2} e^{-\frac{1}{2}\frac{1+r}{1-r}(x^2+y^2)}\\
\times& J_\al(2ixy r^{1/2}(1-r)^{-1}) 
\end{align*} 
(see e.g. (11) of \cite{E97}), 
we can obtain the heat kernel formula for $L_\al$, $\al\ge -1/2$ 
\begin{align}
&e^{-tL_\al}(x,y)=\sum_{k=0}^\iy 
e^{-t(4k+2\al+2)}M^\al_k(x)M^\al_k(y)\notag\\
=&(-i)^\al(\sh 2t)^{-1}e^{-\frac{x^2+y^2}{2}\coth 2t}(xy)^{1/2}
 J_\al\left(\frac{ixy}{\sh 2t}\right), \label{e:etL-ker}
\end{align}
where $J_\al$ is the Bessel function of order $\al$. 
Now if $\al\ge 1/2$, (\ref{e:etL-gb}) follows from the usual asymptotics of $J_\al$ 
(see \cite{Wa45} or \cite[VIII.5]
{St93}):
\begin{align*}
&J_\al(z)=\begin{cases}
z^\al&z\to 0\\
z^{-1/2}&z\to \iy,
\end{cases}\\
&J'_\al(z)=\begin{cases}
z^{\al-1}&z\to 0\\
z^{-1/2}&z\to \iy.
\end{cases}\\
\end{align*}
\end{proof}

\subsection{Laguerre operator in $\R_+^n$} 
A natural extension of $L_\al$ to $n$-dimensions is for $\al=(\al_1,\dots,\al_n)$, $\al_j\ge -\frac12$ 
\begin{align*}
\mcal{L}_\al=\sum_{j=1}^n
(-\pa^2_{x_j})+x_j^2+(\al_j^2-\frac{1}{4})x_j^{-2}
\end{align*} 
in $L^2(\R^n_+)$. 
It can be defined as the Friedrich's extension of the form 
\[ (\mcal{L}_\al f,f)=\sum_{j=1}^n\int_{\R^n_+}
|\pa_{x_j}f(x)|^2+\big(x_j^2+(\al_j^2-\frac{1}{4})x_j^{-2}\big)
|f(x)|^2dx
\]
for $f\in C^\iy_0(\R^n_+)$; hence it's domain is a 
subspace of the Sobolev space $
W^1_2(\R^n)$. 

The ONB in $L^2(\R_+^n)$ consists of 
$\{M_k^\al(x)=M_{k_1}^{\al_1}(x_1)\otimes\cdots\otimes 
M_{k_n}^{\al_n}(x_n), x=(x_1,\dots,x_n)\in\R^n_+, 
k=(k_1,\dots,k_n)\in \N_0^n\}$. 
Hence the heat kernel of $\mcal{L}_\al$ is given by 
\begin{align*}
&e^{-t\mcal{L}_\al}(x,y)=\sum_{k=(k_1,\dots,k_n)}
e^{-t\lam_k^\al}M^\al_k(x)M^\al_k(y)\\
=&(-i)^{|\al|}(\sh 2t)^{-n}
e^{-\frac{|x|^2+|y|^2}{2}\coth 2t}
 \prod_{j=1}^n(x_jy_j)^{1/2}
J_{\al_j}\left(\frac{ix_jy_j}{\sh 2t}\right),
\end{align*}
where $\lam_k^\al=
\sum_{j=1}^n(4k_j+2\al_j+2)$, $x,y\in \R^n_+$.
Therefore applying the 1D result we easily obtain the estimates in (\ref{e:etL-gb}) for
$\mcal{L}_\al$ provided $\al_j\ge 1/2$, $j=1,\dots,n$.



\vs{.05077in}
\nd
\rk Thangavelu \cite{Th06} obtained a heat kernel formula in the special case $\al\in\N_0$ 
via different method. 
We observe that 
the formula in  \cite[Theorem 2.17]{Th06} for 
$L(\al)=-d^2/dx^2{-}(2\al+1)x\inv d/dx+x^{2}/4$ 
acting in the weighted space $(L^2(\R_+), x^{2\al+1})$ is equivalent with 
(\ref{e:etL-ker}) for $L_\al$ acting in the unweighted space $L^2(\R_+)$. 

\subsection{Schr\"odinger operator and associated heat kernel}
Proposition \ref{pr:etH-phi} shows that the upper Gaussian bound estimates of 
$e^{-tH}(x,y)$ imply the decay estimates (\ref{e:der-ker-phi})  in Assumption \ref{a:phi-dec}. 
However for bounded $V\neq 0$ the gradient estimates for $\nabla_xp_t$ are not valid in general, 
this is one reason we work 
in a more direct way to deal with the decay of  $\nabla_x\vphi_j(H)(x,y)$ for all $j\in\Z$ 
as illustrated in the one dimensional case. 
 The 
heat kernel approach seems to work more efficiently for unbounded potentials. 
For bounded potentials, when the Gaussian bounds are not available, we can consider, for instance, the radial case in three dimensions using Volterra type equation
for the eigenfunction of $H$, or 
 the non-radial case 
using stationary phase method \cite{Sch05b,Cu00,DP05,Z07}. 

\section{Conclusion}
The Littlewood-Paley theory of $\dot{B}_p^{\al,q}(H)$ and $\dot{F}_p^{\al,q}(H)$ surveyed in this paper leaves open 
other problems in
the area of harmonic analysis and PDE.  Here we would like to mention a few problems related to our 
subject. 

a. The further characterization and identification
of $\dot{B}_p^{\al,q}(H), \dot{F}_p^{\al,q}(H)$, $0<p,q\le\iy$ as well as other elements of fundamental theory of
function spaces, including distribution theory, have yet to 
be understood, compared with classical Fourier analysis \cite{Tr83,Tr92, FJW}.
It may involve semigroup method, Riesz transform as well as singular integrals for certain class of rough potentials.

b. Applying the Littlewood-Paley decomposition to establish dispersion and Strichartz estimates for the 
perturbed wave, Klein-Gordon and Schr\"odinger equations with potentials.  Although there has been quite extensive work
in this area, cf. \cite{
Be94,Ca03,Cu00,
Sch05c,Y95,
Y05},  regularities  
of these estimates involving the associated  
function spaces 
have not received comparable attention. 
 We hope the development of the Littlewood-Paley theory in our continuing investigation
could give a systematic treatment of the regularity problems, which are related 
to  one of the open problems 
concerning Strichartz type estimates for wave equations with potentials 
in the presence of resonance or eigenvalue at zero energy, cf. \cite{KS07,Sch05b,OOZ06}.

c. For the Laplacian-Beltrami operator  $-\De_g$ on a (complete) Riemannian manifold 
$(M,g)$, in many cases the 
generic heat kernel estimates 
are valid \cite{Da89a, Da89b,Gr95}.  
However for a Schr\"odinger operator $H_V=-\De_g+V$ on $M$,  
the heat kernel estimates 
is in general not valid,  
this makes it unique and more difficult especially for the high and low energy analysis of $H_V$. 
 The treatment here might throw a light on considering the analogous problems in a, b 
on $M$; in particular we are interested in obtaining the analogous results 
 for $H_V$ 
 on Riemannian symmetric spaces. 
We refer to \cite{AO03,
OS92,Sk98b,Sk96,Tat01,BaCS06,
Pie06}  for some of the recent development in this direction.


\begin{thebibliography}{99}
\renewcommand{\baselinestretch}{1}








\bibitem{AO03}
J.-P. Anker, P.~Ostellari, 
The heat kernel on noncompact symmetric spaces. {\em Lie Groups and Symmetric spaces}. 
AMS Translations, 
Ser. 2, volume {\bf 210} (2003), 27-46. 




\bibitem{BaCS06} 
V.~Banica, R.~Carles, G.~Staffilani,  
 Scattering theory for radial nonlinear Schr\"odinger equations on hyperbolic space. http://lanl.arXiv.org/math.AP/0607186.



\bibitem{Be94} M.~Beals, Optimal $L^\infty$ decay for solutions to the wave equation with a potential.  {\em Comm. P.D.E.} {\bf 19} (1994), no. 7-8, 1319-1369.



\bibitem{BZ} 
J.~Benedetto, S.~Zheng, Besov spaces for the 
Schr\"odinger operator with barrier potential.  Submitted. 
{http://lanl.arXiv.org/math.CA/0411348}. 



 
\bibitem{BT06} 
B.~Bongioanni, J.~Torrea,  
Sobolev spaces associated to the harmonic oscillator. 
{\em Proc. Indian Acad. Sci. Math. Sci}. {\bf 116} (2006), no. 3, 337-360. 



                                     

 
 


\bibitem{Ca03}
R.~Carles, Semi-classical Schr\"odinger equations with harmonic potential and nonlinear perturbation. 
{\em Annales I.H.P., Analyse non lin\'eaire} {\bf 20} (2003), no. 3, 501-542. 














\bibitem{Cu00}
S.~Cuccagna, On the wave equation with a potential. {\em Comm. P. D. E.}
{\bf 25} (2000), 1549-1565.



\bibitem{DF06}
P.~D'Ancona, L.~Fanelli,  $L\sp p$-boundedness of the wave operator for the one dimensional Schr\"odinger operator. {\em Comm. Math. Phys}. {\bf 268} (2006), no. 2, 415-438. 



\bibitem{DP05}
P.~D'Ancona, V.~Pierfelice, 
On the wave equation with a large rough potential. {\em J. Funct. Anal}. 
{\bf 227} (2005), no. 1, 30-77.


\bibitem{Da89a}
E.~Davies, {\em Heat Kernels and Spectral Theory}. Cambridge University Press, Cambridge, 1989. 

\bibitem{Da89b} \bysame,
Pointwise bounds on the space and time derivatives of heat kernels. {\em J. Operator Theory}
{\bf 21} (1989), 367-378.





\bibitem{DT79}
P.~Deift, E.~Trubowitz,  Inverse scattering on the line. \emph{Comm. Pure Appl. Math}. vol. XXXII. (1979),  121-251.



\bibitem{DOS02} 
X.~Duong, E.~Ouhabaz, A.~Sikora, Plancherel type
estimates and sharp spectral multipliers. {\em J. Funct. Anal}. {\bf 196} (2002).  443-485.

\bibitem{DY05}
X.~Duong, L.~Yan, 
Duality of Hardy and BMO spaces associated with operators with heat kernel bounds.
{\em J. Amer. Math. Soc}. {\bf 18} (2005), 943-973. 

\bibitem{D97} J.~Dziuba\'nsk, Triebel-Lizorkin spaces associated with 
Laguerre and Hermite expansions. {\em Proc. Amer. Math. Soc.}{\bf
125} (1997), 3547-3554.



\bibitem{D01} \bysame,  
A spectral multiplier theorem for $H\sp 1$ spaces associated with Schr\"odinger operators with potentials satisfying a reverse H\"older inequality. {\em Illinois J. Math}. {\bf 45} (2001), no. 4, 1301-1313. 

\bibitem{DZ05}
J.~Dziuba\'nski, J.~Zienkiewicz, 
Hardy spaces $H\sp 1$ for Schr\"odinger operators with compactly supported potentials. 
{\em Ann. Mat. Pura Appl}. (4) {\bf 184} (2005), no. 3, 315-326.



\bibitem{DZ02}\bysame,
$H^p$  spaces for Schr\"odinger operators, in {\em Fourier Analysis and Related Topics}, 
Banach Center Publ. {\bf 56}, Inst. Math., Polish Acad. Sci., 2002, 45-53.


\bibitem{E95} J.~Epperson, Triebel-Lizorkin Spaces for Hermite
expansions.  {\em Studia Math.} {\bf 114} (1995), no.1, 87-103.

\bibitem{E96}\bysame, Hermite multipliers and pseudo-multipliers. {\em
Proc. Amer. Math. Soc.} {\bf 124} (1996), no. 7, 2061-2068.

\bibitem{E97}\bysame, 
Hermite and Laguerre wave packet expansions. {\em Studia
Math.} {\bf 126} (1997), no. 3, 199-217.










\bibitem{FJW} 
M.~Frazier, B.~Jawerth, G.~Weiss, {\em Littlewood-Paley Theory
and the Study of Function Spaces},  Conference Board of the
Math. Sci. {\bf 79}, 1991.






\bibitem{GS04} 
 M.~Goldberg, W.~Schlag, Dispersive estimates for Schr\"odinger operators in dimensions one and three. \emph{Comm. Math. Phys}. {\bf 251} (2004), no. 1, 157-178.
 


\bibitem{Gr95}
A.~Grigor'yan, Upper bounds of derivatives of the heat kernel on an arbitrary complete manifold. 
{\em J. Funct. Anal}. {\bf 127} (1995), no. 2, 363--389. 




\bibitem{GV95} 
J.~Ginibre, G.~Velo, Generalized Strichartz inequalities for
the wave equation, {\em J. Funct. Anal}. {\bf 133} (1995), 50-68.





\bibitem{He90a} 
W.~Hebisch: A multiplier theorem for Schr\"odinger operators.
\textit{Colloq. Math}, {\bf 60/61} (1990),  no. 2, 659-664.


\bibitem{He95}\bysame,  
Functional calculus for slowly decaying kernels. Preprint, 1995.\\ 
http://www.math.uni.wroc.pl/\symbol{126}hebisch








\bibitem{Hor60} 
L.~H\"ormander,  Estimates for translation invariant operators in
$L^p$ spaces. {\em Acta Math}. {\bf 104} (1960), 93-140.





\bibitem{JN94} 
A.~Jensen, S.~Nakamura, Mapping properties of functions
of Schr\"odinger operators between $L^p$ spaces and Besov spaces. 
 {\em Spectral and Scattering Theory and Applications}, Advanced 
Studies in Pure Math. {\bf 23} (1994), 187-209.









\bibitem{KT98}
M.~Keel, T.~Tao, Endpoint Strichartz estimates. 
{\em Amer. J. Math}. {\bf 120} (1998), 955-980. 






\bibitem{KS07}
J.~Krieger, W.~Schlag, On the focusing critical semi-linear wave equation. {\em Amer. J. Math}. {\bf 129} (2007), no. 3, 843-913. 




 
 
















\bibitem{OS92}
G.~\'Olafsson, H. Schlichtkrull, Wave propagation on Riemannian symmetric 
space. {\em J. Funct. Anal}. {\bf 107} (1992), 270-278.

\bibitem{OZ06}
G.~\'Olafsson, S.~Zheng, Function spaces associated with 
Schr\"odinger operators: the P\"oschl-Teller potential.  
{\em Jour. Four. Anal. Appl}. {\bf 12} (2006), no.6, 653-674. 

\bibitem{OOZ06}
G.~\'Olafsson, K.~Oskolkov, S.~Zheng, Spectral multipliers for  
Schr\"odinger operators II. Preprint.








\bibitem{Pie06} 
V.~Pierfelice, Weighted Strichartz estimates for the radial 
perturbed Schr\"odinger equation on the hyperbolic space. 
{\em Manuscripta Math}. {\bf 120}, (2006), 377-389.







\bibitem{Sch05b} 
W.~Schlag, A remark on Littlewood-Paley theory for the distorted 
Fourier transform. 
{\em Proc. Amer. Math. Soc}. {\bf 135} (2007), no. 2, 437-451.

\bibitem{Sch05c} 
\bysame, Spectral theory and nonlinear partial differential equations: a survey. 
{\em Discrete Contin. Dyn. Syst}. {\bf 15} (2006), no. 3, 703-72. 
{\em http://lanl.arXiv.org/math.AP/0509019}.












\bibitem{Si82}
B.~Simon, Schr\"odinger semigroups. {\em Bull. Amer. Math. Soc}. {\bf 7}
(1982), no.3, 447-526.




\bibitem{Sk98b} L.~Skrzypczak, 
The Triebel-Lizorkin scale of function spaces for the Fourier-Helgason transform.  {\em Math. Nachr}.  {\bf 190}  (1998), 251-274. 

\bibitem{Sk96}\bysame, 
Heat semi-group and function spaces on symmetric spaces on non-compact type. {\em Z. Anal. Anwendungen} {\bf 15}  (1996),  no. 4, 881--899. 


 





\bibitem{St93}
E.~Stein, {\em Harmonic analysis, Real-variable methods, orthogonality, and 
oscillatory integrals}. Princeton Univ. Press, 1993.

\bibitem{ST02} 
G.~Staffilani, D.~Tataru, 
Strichartz estimates for a Schr\"odinger operator with nonsmooth coefficients. {\em Comm. P. D. E.}  {\bf 27} (2002),  no. 7-8, 1337-1372.

\bibitem{Str77}
R.~Strichartz, Restrictions of Fourier transforms to quadratic surfaces and decay of solutions of
wave equations. {\em Duke Math. J}.  {\bf 44} (1977), no.3, 705-714 .

\bibitem{Sz06a}
Z.~Szabo, De Broglie geometry on Zeeman manifolds: 
A new non-perturbative approach to the infinities of QED. 
{\em The 26th International Conference on Group Theoretical Methods in Physics} (2006). http://lanl.arXiv.org/math-ph/0609027.





\bibitem{Tat01} 
D.~Tataru,  Strichartz estimates in the hyperbolic space and global existence for the semilinear wave equation.
{\em Trans. Amer. Math. Soc}. {\bf 353} (2001), 795-807. 





\bibitem{Th93} S.~Thangavelu, {\em Lectures on Hermite and Laguerre
expansions}. Princeton Univ. Press, 1993.




\bibitem{Th06}\bysame,  
Hermite and Laguerre semigroups: Some recent developments. 
To appear in CIMPA lecture notes, 2006.



\bibitem{Tie06}
J.~Tie, The twisted Laplacian on $\mathbb{C}^n$ 
and the sub-Laplacian on $H^n$. {\em Comm. P. D. E}.  {\bf 31}  (2006),  no. 7-9, 1047-1069. 




\bibitem{Tr83} H.~Triebel, {\em Theory of Function Spaces}. 
Birkh\"{a}user Verlag, 1983.

\bibitem{Tr92}\bysame, 
{\em Theory of Function Spaces II}. Monographs
Math. {\bf 84}, Birkh\"auser, Basel, 1992.







\bibitem{Wa45} G.~Watson, {\em A Treatise on the 
Theory of Bessel Functions}. Cambridge University Press, 
Second Ed. 1944. 




\bibitem{Y95} K.~Yajima, The $W^{k,p}$-continuity of wave operators for Schr\"odinger operators. {\em J. Math. Soc. Japan} {\bf 47} (1995), 551-581.



\bibitem{Y05} \bysame, 
Dispersive estimate for Schr\"odinger equations with threshold
resonance and eigenvalue.  {\em Comm. Math. Phys}. {\bf 259} (2005), 475-509. 


 
\bibitem{Z04} 
S.~Zheng, A representation formula related to Schr\"odinger 
operators.  {\em Anal.~Theory~Appl.} {\bf 20} (2004), no.3, 294-296. 


\bibitem{Z06a} 
\bysame,  Littlewood-Paley theorem for  Schr\"odinger operators. 
{\em Anal. Theory Appl}. {\bf 22} (2006), no.4, 353-361.

\bibitem{Z06b}\bysame, Spectral multipliers for Schr\"odinger operators I.  Preprint. \n
http://lanl.arXiv.org/math.AP/0610096.






\bibitem{Z07}\bysame, 
{Spectral calculus, function spaces and dispersive equation with 
a critical potential}. In preparation.






\end{thebibliography}
\end{document}